\documentclass{amsproc}
\usepackage{amssymb}
\usepackage{amsmath}
\usepackage{eurosym}
\usepackage{amsfonts}

            \advance \topmargin by -\headheight \advance
            \topmargin by -\headsep
            \evensidemargin
            \oddsidemargin
\newtheorem{theorem}{Theorem}[section]
\newtheorem{lemma}[theorem]{Lemma}

\newtheorem{remark}[theorem]{Remark}
\numberwithin{equation}{section}
\input{tcilatex}

\begin{document}
\title[THE BOUNDEDNESS OF CAUCY INTEGRAL OPERATOR ]{THE BOUNDEDNESS OF
CAUCHY INTEGRAL OPERATOR ON A DOMAIN HAVING CLOSED ANALYTIC BOUNDARY}
\thanks{}
\author[Y\"{u}ksel ~Soykan]{Y\"{U}KSEL SOYKAN}
\maketitle

\begin{center}
\textsl{Zonguldak B\"{u}lent Ecevit University, Department of Mathematics, }

\textsl{Art and Science Faculty, 67100, Zonguldak, Turkey }

\textsl{e-mail: \ yuksel\_soykan@hotmail.com}
\end{center}

\textbf{Abstract.} In this paper, we prove that the Cauchy integral
operators (or Cauchy transforms) define continuous linear operators on the
Smirnov classes for some certain domain with closed analytic boundary.



\section{Introduction}

As usual, we define the Hardy space $H^{2}=H^{2}(\Delta )$ as the space of
all functions $f:z\rightarrow \sum_{n=0}^{\infty }a_{n}z^{n}$ for which the
norm $\left( \left\Vert f\right\Vert =\sum_{n=0}^{\infty }\left\vert
a_{n}\right\vert ^{2}\right) ^{1/2}$ is finite. Here, $\Delta $ is the open
unit disc. For a more general simply-connected domain $\Omega $ in the
complex plane $\mathbb{C}$ with at least two boundary points, and a
conformal mapping $\varphi $ from $\Omega $ onto $\Delta $ (that is, a
Riemann mapping function), a function $g$ analytic in $\Omega $ is said to
belong to the Smirnov class $E^{2}(\Omega )$ if and only if $\ g=(f\circ
\varphi )\varphi ^{\prime }{}^{1/2}$ $\ $for some $f\in H^{2}(\Delta )$
where $\varphi ^{\prime }{}^{1/2}$ is an analytic branch of the square root
of $\varphi ^{\prime }.$ The reader is referred to [\ref{bib:dur}], [\ref%
{bib:dur1}], [\ref{bib:go}], [\ref{bib:kh1}], and references therein for a
basic account of the subject. $\partial \Delta $ and $\partial \Omega $ will
be used to denote the boundary of open unit disc $\Delta $ and the boundary
of $\Omega $ respectively.

\smallskip Suppose that $\Gamma $ is a simple $\sigma $-rectifiable arc (not
necessarily closed). The notation $L^{p}(\Gamma )$ will denote the $L^{p}$
space of normalized arc length measure on $\Gamma $. Let $\Omega $ denote
the complement of $\Gamma $. The Cauchy Integral of a function $\widetilde{f}
$ defined on $\Gamma $ and integrable relative to arc length is defined as:%
\begin{equation}
C_{\Omega }\widetilde{f}(z)=\frac{1}{2\pi i}\int_{\Gamma }\frac{\widetilde{f}%
(\zeta )}{\zeta \text{ }-z}d\zeta \text{ \hspace{0in}\qquad }(z\in \Omega ).
\label{equ:yyyy}
\end{equation}%
\smallskip $C_{\Omega }\widetilde{f}$ is analytic at each point of $\Omega $%
%
%
%
%
%
%
%
%
%
%
%
%
%
%
%
%
. If $\Gamma $ is not closed, then (\ref{equ:yyyy}) defines a single
analytic function. If $\Gamma $ is closed, then $\Omega $ has two
components, the interior and the exterior of $\Gamma $. Then in each
component of (\ref{equ:yyyy}) defines an analytic function%
.

Recall that a closed analytic curve is a curve $\gamma =k(\partial \Delta )$
where $k$ is analytic and conformal in a neighbourhood $U$ of $\partial
\Delta $. If $\gamma $ is simple it is called an analytic Jordan curve.

In this paper we prove

\begin{theorem}
\label{the:cauchyint}Suppose that that $D$ is a bounded simply connected
domain and $\gamma =$ $\partial D$ is a closed analytic curve (e.g.
ellipse). Then the Cauchy Integral 
\begin{equation*}
C_{D}\widetilde{f}(z)=\frac{1}{2\pi i}\int_{\partial D}\frac{\widetilde{f}%
(\zeta )}{\zeta \text{ }-z}d\zeta
\end{equation*}%
defines a continuous linear operator mapping $L^{2}(\partial D)$ into $%
E^{2}(D)$.
\end{theorem}

\begin{remark}
The result of the above Theorem is well known in the literature; see, for
example, [\ref{Coifman}] and [\ref{Dyn}]. However, we give a basic and
direct proof.
\end{remark}

To prove Theorem \ref{the:cauchyint} we need the following Lemma and Remark.

\begin{lemma}[{[\protect\ref{bib:dur}], p.170}]
Suppose that $\Omega $ is a simply connected and bounded domain and the
boundary $\Gamma =\partial \Omega $ of $\Omega $ is a rectifiable Jordan
curve. Then

\begin{description}
\item[i.] Each $f\in E^{2}(\Omega )$ has a nontangential limit function $%
\widetilde{f}\in L^{2}(\partial \Omega )$ and 
\begin{equation*}
||f||_{E^{2}(\Omega )}^{2}=||\widetilde{f}||_{L^{2}(\partial \Omega )}^{2}=%
\frac{1}{2\pi }\int_{\partial \Omega }|\widetilde{f}(z)|^{2}|dz|.
\end{equation*}

\item[ii.] Each $f\in E^{2}(\Omega ))$ has a Cauchy representation%
\begin{equation}
f(z)=\frac{1}{2\pi i}\int_{\partial \Omega }\frac{\widetilde{f}(\zeta )}{%
\zeta -z}\Omega \zeta \qquad \qquad (z\in \Omega ).  \label{equ:tutu}
\end{equation}
\end{description}
\end{lemma}

In this case, for equation \ref{equ:tutu}, we say that Cauchy Integral
Formula is valid.

A special case of the above theorem is the following Remark. In fact, we
prove Theorem \ref{the:cauchyint} using this Remark.

\begin{remark}[{[\protect\ref{bib:li2}], p.423}]
The Cauchy integral formula 
\begin{equation*}
C_{\Delta }\widetilde{f}(z)=\frac{1}{2\pi i}\int_{\partial \Delta }\frac{%
\widetilde{f}(\zeta )}{\zeta -z}d\zeta
\end{equation*}%
defines a continuous linear operator $C_{\Delta }:L^{2}(\partial \Delta
)\rightarrow E^{2}(\Delta )$ with $||C_{\Delta }||=1.$
\end{remark}

There is another integral operator $C_{\Gamma }$ on the curve $\Gamma ,$
which is sometimes called also Cauchy integral, viewed as an operator-valued
function of the curve. This operator is given by a principal value singular
integral: if $\widetilde{f}$ is a function on $\Gamma $, we define $%
C_{\Gamma }(\widetilde{f})$ on $\Gamma $ by%
\begin{equation*}
C_{\Gamma }\widetilde{f}(z)=\frac{1}{2\pi i}P.V.\int_{\Gamma }\frac{%
\widetilde{f}(\zeta )}{\zeta \text{ }-z}d\zeta \text{ \hspace{0in}\qquad }%
(z\in \Gamma ).
\end{equation*}%
The operator $C_{\Gamma }$ is probably less familiar than $C_{\Omega }.$
These operators are very important in real and complex analysis, and have
attracted many mathematicians to investigate them.

In fact, there are other types of Cauchy integrals and there have been
extensive literature about them and theirs applications as papers and books.
For the books concerning Cauchy type integrals and related subjects; see,
for instance, [\ref{bib:Bell}], [\ref{bib:Cima}], [\ref{bib:Murai}], [\ref%
{bib:Tolsa1}] and [\ref{bib:Tolsa2}]. For the books concerning boundedness
of Cauchy type integrals; see, for example, [\ref{bib:Edmunds}], [\ref%
{bib:Halmos}] and [\ref{bib:Okikiolu}].

\section{\protect\smallskip Proof of Theorem \protect\ref{the:cauchyint}}

We are now ready to give the proof of Theorem \ref{the:cauchyint}.

\proof
Suppose that $\varphi $ is a conformal map of $D$ onto $\Delta $. Let $\psi
=\varphi ^{-1}:\Delta \rightarrow D$. \smallskip \smallskip \smallskip
\smallskip Consider the maps $C_{\Delta }:L^{2}(\partial \Delta )\rightarrow
E^{2}(\Delta )$ is given by $C_{\Delta }\widetilde{f}=\frac{1}{2\pi i}%
\int_{\partial \Delta }\frac{\widetilde{f}(\zeta )}{\zeta -z}d\zeta $ and $%
\widetilde{U}_{\psi }$\smallskip :$L^{2}(\partial D)\rightarrow
L^{2}(\partial \Delta )$ is given by $\widetilde{U}_{\psi }f(z)=f(\psi
(z))\psi ^{\prime }(z)^{1/2}$ and $U_{\varphi }$\smallskip :$E^{2}(\Delta
)\rightarrow E^{2}(D)$ is given by $U_{\varphi }f(z)=f(\varphi (z))\varphi
^{\prime }(z)^{1/2}$. \smallskip $U_{\varphi }$ and $\widetilde{U}_{\psi }$
are unitary operators. The situation is illustrated in the Figure \ref%
{fig:aabbccddf}.

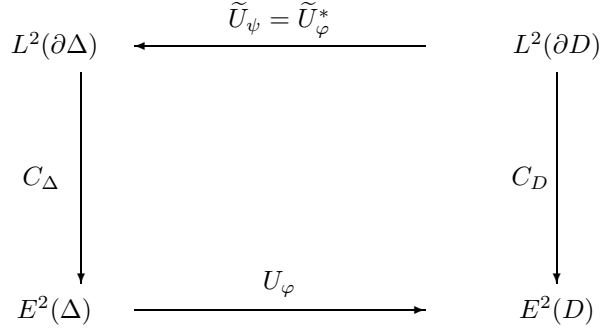
\begin{figure}[h]
\begin{picture}(230,170)(-70,17)
\put(30,130){\makebox(0,0){$L^2(\partial \Delta )$}}
\put(170,130){\vector(-1,0){110}}
\put(115,140){\makebox(0,0){$\widetilde{U}_{\psi} ={\widetilde{U}_\varphi }^{*}$}}
\put(220,130){\makebox(0,0){$L^2(\partial D)$}}
\put(210,80){\makebox(0,0){$C_D$}}
\put(220,120){\vector(0,-1){80}}
\put(220,30){\makebox(0,0){$E^2(D)$}}
\put(115,40){\makebox(0,0){$U_{\varphi} $}}
\put(60,30){\vector(1,0){110}}
\put(40,120){\vector(0,-1){80}}
\put(25,80){\makebox(0,0){$C_\Delta $}}
\put(30,30){\makebox(0,0){$E^2(\Delta )$}}

\end{picture}
\newline
\caption{the maps $U_{\protect\varphi }$ and $\protect\widetilde{U}_{\protect%
\psi }$ }
\label{fig:aabbccddf}
\end{figure}
For $\widetilde{f}\in L^{2}(\partial \Delta ),$ we have $U_{\psi }C_{D}%
\widetilde{U}_{\varphi }\widetilde{f}=U_{\psi }C_{D}(\widetilde{f}\circ
\varphi )\varphi ^{\prime ^{1/2}}$ and%
\begin{equation*}
C_{D}(\widetilde{f}\circ \varphi (w))\varphi ^{\prime }(w)^{1/2}=\frac{1}{%
2\pi i}\int_{\partial D}\frac{\widetilde{f}(\varphi (\zeta ))\varphi
^{\prime }(\zeta )^{1/2}}{\zeta -w}d\zeta
\end{equation*}%
so that%
\begin{eqnarray*}
U_{\psi }C_{D}\widetilde{U}_{\varphi }\widetilde{f}(z) &=&\psi ^{\prime
}(z)^{1/2}\frac{1}{2\pi i}\int_{\partial D}\frac{\widetilde{f}(\varphi
(\zeta ))\varphi ^{\prime }(\zeta )^{1/2}}{\zeta -\psi (z)}d\zeta \\
&=&\frac{1}{2\pi i}\int_{\partial \Delta }\frac{\widetilde{f}(w)\psi
^{\prime }(w)^{1/2}\psi ^{\prime }(z)^{1/2}}{\psi (w)-\psi (z)}dw\qquad
(z\in \Delta ).
\end{eqnarray*}

Then we have kernel%
\begin{equation*}
K(z,w)=\frac{1}{w-z}H(z,w)
\end{equation*}%
where 
\begin{equation*}
H(z,w)=\frac{(w-z)\psi ^{\prime }(w)^{1/2}\psi ^{\prime }(z)^{1/2}}{\psi
(w)-\psi (z)}.
\end{equation*}

\smallskip 
For any $1<r$, denote $\Delta _{r}$ and $\sigma _{r}$ by $\Delta
_{r}=\{z:\left\vert z\right\vert <1\}$ and $\sigma _{r}=\partial \Delta
_{r}=\{z:\left\vert z\right\vert =r\}$.

Since $\partial D$ is analytic (closed) curve, $\psi $ is analytic and
conformal in a neighbourhood of $\overline{\Delta }.$ So without loss of
generality, we may assume that $\psi $ is analytic and conformal on $\Delta
_{R}$ for some $R>1.$ Then $\psi $ is analytic and conformal on and inside \-%
$\sigma _{r^{\prime }}$ where $R>r^{\prime }>1$.

\smallskip Choose $r,s$ such that $1<r<s<R.$ We shall show that $H$ is
analytic on $\Delta _{r}\times \Delta _{r}$. Hence, because $r$ is
arbitrary, it will follow that $H$ is analytic on $\Delta _{R}\times \Delta
_{R}$.

Fix $z\in $ $\Delta _{R}$. Then $w\rightarrow F_{z}(w)=H(z,w)$ is analytic
in $\Delta _{R}$ except at $w=z$. But since residue at $w=z$ is 
$0$, the singularities of $H$ for $w=z$ is removable. %
Hence $w\rightarrow H(z,w)$ is analytic on $\Delta _{R}\supseteq \overline{%
\Delta _{s}}\supseteq \overline{\Delta _{r}}$ 
. We can thus apply Cauchy's integral formula to it, giving 
\begin{equation}
H(z,w)=\frac{1}{2\pi i}\int_{\sigma _{r}}H(z,v)\frac{1}{v-w}dv\text{ \ }(%
\text{for }w\in \Delta _{r}\text{ and fixed }z\in \Delta _{s},v\in \sigma
_{r}).  \label{equ:ttyy}
\end{equation}%
Hence since $z\in \Delta _{s}$ is arbitrary, for all $z\in \Delta _{s}$ and $%
w\in \Delta _{r}$, equation \ref{equ:ttyy} is valid.

\smallskip By symmetry
, for every $v\in $ $\Delta _{R}$ the function $z\rightarrow H(z,v)$ is
analytic on $\Delta _{R}$. Hence%
\begin{equation}
H(z,v)=\frac{1}{2\pi i}\int_{\sigma _{s}}H(u,v)\frac{1}{u-z}du\qquad (z\in
\Delta _{s},v\in \Delta _{R},u\in \sigma _{s}).  \label{equ:tryt}
\end{equation}

$H$ is separately continuous on $\Delta _{R}\times \Delta _{R}$, so it is
(jointly) continuous on $\Delta _{R}\times \Delta _{R}$. Substitute the
value of $H(z,v)$ from (\ref{equ:tryt}) in the integrand of (\ref{equ:ttyy}%
). Since the function $H(u,v)$ is continuous, we obtain 
\begin{equation*}
H(z,w)=\left( \frac{1}{2\pi i}\right) ^{2}\int_{\sigma _{r}}\int_{\sigma
_{s}}\frac{H(u,v)dudv}{(u-z)(v-w)}\qquad (z\in \Delta _{s},w\in \Delta _{r},%
\text{ }v\in \sigma _{r},u\in \sigma _{s})
\end{equation*}%
Since $r$ is arbitrary, we will show that $H(z,w)$ is analytic on $\Delta
_{R}\times $ $\Delta _{R}.$

Now,%
\begin{equation*}
\frac{1}{(u-z)(v-w)}=\sum_{m,n=0}^{\infty }\frac{z^{m}}{u^{m+1}}\frac{w^{n}}{%
v^{n+1}}\qquad (v\in \sigma _{r},u\in \sigma _{s})
\end{equation*}%
and 
this series is uniformly convergent for $z\in \Delta _{s},w\in \Delta _{r},$ 
$v\in \sigma _{r},u\in \sigma _{s}$ 
.

\smallskip Hence%
\begin{equation}
H(z,w)=\left( \frac{1}{2\pi i}\right) ^{2}\int_{\sigma _{r}}\int_{\sigma
_{s}}\sum_{m,n=0}^{\infty }\frac{z^{m}}{u^{m+1}}\frac{w^{n}}{v^{n+1}}%
H(u,v)dudv\qquad (z\in \Delta _{s},w\in \Delta _{r},\text{ }v\in \sigma
_{r},u\in \sigma _{s})  \label{equ:ttss}
\end{equation}

Since $H$ is bounded , the series $\sum_{m,n=0}^\infty \frac{z^m}{u^{m+1}}%
\frac{w^n}{v^{n+1}}H(u,v)$ is uniformly convergent for $z\in \Delta _s,w\in
\Delta _r,$ $v\in \sigma _r,u\in \sigma _s$. %

Because of uniformly convergence%
, we can integrate the series (\ref{equ:ttss}) term-by-term and we obtain 
\begin{eqnarray}
H(z,w) &=&\sum_{m,n=0}^{\infty }z^{m}w^{n}\left( \frac{1}{2\pi i}\right)
^{2}\int_{\sigma _{r}}\int_{\sigma _{s}}\frac{H(u,v)}{u^{m+1}v^{n+1}}%
dudv\qquad (z\in \Delta _{s},w\in \Delta _{r},\text{ }v\in \sigma _{r},u\in
\sigma _{s})  \notag \\
&=&\sum_{m,n=0}^{\infty }a_{mn}z^{m}w^{n}  \label{equ:tyuy}
\end{eqnarray}%
where the coefficients $a_{mn}$ are given by the integral formula%
\begin{equation*}
a_{mn}=\left( \frac{1}{2\pi i}\right) ^{2}\int_{\sigma _{r}}\int_{\sigma
_{s}}\frac{H(u,v)}{u^{m+1}v^{n+1}}dudv.
\end{equation*}

%
%
%
%
\smallskip Since $H$ is bounded, we obtain 
\begin{eqnarray*}
\left\vert a_{mn}\right\vert &=&\left\vert \left( \frac{1}{2\pi i}\right)
^{2}\int_{\sigma _{r}}\int_{\sigma _{s}}\frac{H(u,v)}{u^{n+1}v^{m+1}}%
dudv\right\vert \\
&\leq &\frac{1}{4\pi ^{2}}(2\pi s)(2\pi r)\left\Vert H\right\Vert _{\infty }%
\frac{1}{s^{m+1}r^{n+1}} \\
&\leq &s^{2}\left\Vert H\right\Vert _{\infty }\frac{1}{r^{m+1}r^{n+1}}
\end{eqnarray*}%
where $\left\Vert H\right\Vert _{\infty }=\sup_{u,v\in \sigma ^{\prime
}}\left\vert H(u,v)\right\vert <\infty $ and 
so%
\begin{equation*}
\sum_{m,n=0}^{\infty }\left\vert a_{mn}\right\vert <\infty
\end{equation*}%
(i.e. the series $\sum_{m,n=0}^{\infty }a_{mn}$ is absolutely convergent ).
Hence $\sum_{m,n=0}^{\infty }a_{mn}z^{m}w^{n}$ is absolutely convergent on $%
\Delta _{r}\times \Delta _{r}$. Thus $H(z,w)$ is analytic on $\Delta
_{r}(\subseteq \Delta _{s})\times \Delta _{r}$ and so since $r$ is arbitrary
it is analytic on $\Delta _{R}\times \Delta _{R}.$

Now the series $H(z,w)=\sum_{m,n=0}^{\infty }a_{mn}z^{m}w^{n}$ is uniformly
convergent for $z\in \Delta _{s},w\in \Delta _{r}$. If we set $A=U_{\psi
}C_{D}\widetilde{U}_{\varphi }$, then we have%
\begin{eqnarray}
Af(z) &=&U_{\psi }C_{D}\widetilde{U}_{\varphi }f(z)\text{ }\qquad (f\in
L^{2}(\partial \Delta ),z\in \Delta )  \notag \\
&=&\frac{1}{2\pi i}\int_{\partial \Delta }\frac{f(w)\psi ^{\prime
}(w)^{1/2}\psi ^{\prime }(z)^{1/2}}{\psi (w)-\psi (z)}dw\qquad \qquad (w\in
\partial \Delta )  \notag \\
&=&\frac{1}{2\pi i}\int_{\partial \Delta }\frac{1}{w-z}H(z,w)f(w)dw  \notag
\\
&=&\frac{1}{2\pi i}\int_{\partial \Delta }\frac{1}{w-z}\sum_{m,n=0}^{\infty
}a_{mn}z^{m}w^{n}f(w)dw.  \label{equ:rtytr}
\end{eqnarray}

\smallskip In fact, we will show that the sum and the integral in the
equation \ref{equ:rtytr} can be permutable. For $f\in L^{2}(\partial \Delta
),z\in \Delta $, we have 
\begin{equation*}
\sum_{m,n=0}^{\infty }\int_{\partial \Delta }\left\vert \frac{1}{w-z}%
a_{mn}z^{m}w^{n}f(w)\frac{dw}{2\pi i}\right\vert \leq \sum_{m,n=0}^{\infty
}\int_{\partial \Delta }\frac{1}{\left\vert w-z\right\vert }\left\vert
a_{mn}\right\vert \left\vert z\right\vert ^{m}\left\vert w\right\vert
^{n}\left\vert f(w)\right\vert \frac{\left\vert dw\right\vert }{2\pi }
\end{equation*}%
\begin{eqnarray*}
&\leq &\sum_{m,n=0}^{\infty }\left\vert a_{mn}\right\vert \left\vert
z\right\vert ^{m}\left( \int_{\partial \Delta }\left\vert f(w)\right\vert
^{2}\frac{\left\vert dw\right\vert }{2\pi }\right) ^{1/2}\left(
\int_{\partial \Delta }\frac{1}{\left\vert w-z\right\vert ^{2}}\frac{%
\left\vert dw\right\vert }{2\pi }\right) ^{1/2} \\
&\leq &\sum_{m,n=0}^{\infty }\left\vert a_{mn}\right\vert \left\vert
z\right\vert ^{m}\left\Vert f\right\Vert \left( \frac{1}{(1-\left\vert
z\right\vert )^{2}}\right) ^{1/2}<\infty
\end{eqnarray*}%
and by Tonelli Theorem
\begin{equation}
\frac{1}{2\pi i}\int_{\partial \Delta }\frac{1}{w-z}\sum_{m,n=0}^{\infty
}a_{mn}z^{m}w^{n}f(w)dw=\sum_{m,n=0}^{\infty }a_{mn}z^{m}\frac{1}{2\pi i}%
\int_{\partial \Delta }\frac{1}{w-z}w^{n}f(w)dw  \label{equ:yyttr}
\end{equation}%
and so%
\begin{equation*}
Af(z)=\sum_{m,n=0}^{\infty }a_{mn}z^{m}\frac{1}{2\pi i}\int_{\partial \Delta
}\frac{1}{w-z}w^{n}f(w)dw\qquad (f\in L^{2}(\partial \Delta ),z\in \Delta ).
\end{equation*}

Now consider the following series of operators%
\begin{equation*}
\sum_{m,n=0}^{\infty }a_{mn}M_{m}C_{\Delta }N_{n}
\end{equation*}%
where $M_{m}:H^{2}(\Delta )\rightarrow H^{2}(\Delta )$ is defined by $%
M_{m}f(z)=z^{m}f(z)$ so that $\left\Vert M_{m}\right\Vert =\left\Vert
z^{m}\right\Vert _{\infty }=1$, and $N_{n}:L^{2}(\partial \Delta
)\rightarrow L^{2}(\partial \Delta )$ is defined by $N_{n}f(z)=z^{n}f(z)$ so
that $\left\Vert N_{n}\right\Vert =\left\Vert z^{n}\right\Vert _{\infty }=1$%
%
%
%
%
%
%
%
%
%
%
%
%
%
%
%
%
. So we have%
\begin{equation*}
Af(z)=\sum_{m,n=0}^{\infty }a_{mn}M_{m}C_{\Delta }N_{n}f(z)\qquad (f\in
L^{2}(\partial \Delta ),z\in \Delta )
\end{equation*}

Then the series $A=\sum_{m,n=0}^{\infty }a_{mn}M_{m}C_{\Delta }N_{n}$ is
absolutely convergent in operator norm in the space $B(L^{2}(\partial \Delta
),H^{2}(\Delta ))$, in fact, 
\begin{eqnarray*}
\left\Vert M_{m}C_{\Delta }N_{n}f\right\Vert _{H^{2}(\Delta )} &\leq
&\left\Vert M_{m}\right\Vert \left\Vert C_{\Delta }\right\Vert \left\Vert
N_{n}\right\Vert \left\Vert f\right\Vert \qquad (f\in L^{2}(\partial \Delta
)) \\
&=&\left\Vert f\right\Vert _{L^{2}(\partial \Delta )}
\end{eqnarray*}%
and%
\begin{equation*}
\left\Vert M_{m}C_{\Delta }N_{n}\right\Vert \leq 1
\end{equation*}%
so that 
\begin{equation*}
\sum_{m,n=0}^{\infty }\left\Vert a_{mn}M_{m}C_{\Delta }N_{n}\right\Vert \leq
\sum_{m,n=0}^{\infty }\left\vert a_{mn}\right\vert <\infty
\end{equation*}%
(i.e. $\sum_{m,n=0}^{\infty }a_{mn}M_{m}C_{\Delta }N_{n}$ converges
absolutely) and 
\begin{eqnarray*}
\left\Vert Af\right\Vert _{E^{2}(\Delta )} &=&\left\Vert
\sum_{m,n=0}^{\infty }a_{mn}M_{m}C_{\Delta }N_{n}f\right\Vert _{E^{2}(\Delta
)} \\
&\leq &\sum_{m,n=0}^{\infty }\left\Vert a_{mn}M_{m}C_{\Delta
}N_{n}f\right\Vert \text{ \ }(\text{since }\sum_{m,n=0}^{\infty
}a_{mn}M_{m}C_{\Delta }N_{n}\text{ converges absolutely}) \\
&=&\sum_{m,n=0}^{\infty }\left\vert a_{mn}\right\vert \left\Vert
M_{m}\right\Vert \left\Vert C_{\Delta }\right\Vert \left\Vert
N_{n}\right\Vert \left\Vert f\right\Vert \leq \sum_{m,n=0}^{\infty
}\left\vert a_{mn}\right\vert .1.1.1\left\Vert f\right\Vert \\
&\leq &\left\Vert f\right\Vert \sum_{m,n=0}^{\infty }\left\vert
a_{mn}\right\vert <\infty
\end{eqnarray*}%
and so 
\begin{equation*}
\left\Vert A\right\Vert =\left\Vert \sum_{m,n=0}^{\infty
}a_{mn}M_{m}C_{\Delta }N_{n}\right\Vert \leq \sum_{m,n=0}^{\infty
}\left\vert a_{mn}\right\vert .
\end{equation*}%
This shows that $A=U_{\psi }C_{D}\widetilde{U}_{\varphi }$ is a continuous
operator. It follows that $C_{D}$ is a continuous operator.

\textbf{Second proof of of the continuity of }$A$

\smallskip Since in a Banach space $X$ (here $X=B(L^{2}(\partial \Delta
),H^{2}(\Delta ))$ ), every absolutely convergent series is convergent, in
the norm of $X$, to an element of $X$, 
$\sum_{m,n=0}^{\infty }a_{mn}M_{m}C_{\Delta }N_{n}$ converges to an element $%
B$ $\in B(L^{2}(\partial \Delta ),H^{2}(\Delta ))$, in the sense that 
\begin{equation*}
\lim_{m,n\rightarrow \infty
}\sum_{k=0}^{m}\sum_{l=0}^{n}a_{kl}M_{k}C_{\Delta }N_{l}=B
\end{equation*}%
i.e. 
\begin{equation*}
\left\Vert B-\sum_{k=0}^{m}\sum_{l=0}^{n}a_{kl}M_{k}C_{\Delta
}N_{l}\right\Vert \rightarrow 0\quad \text{as }m,n\rightarrow \infty .
\end{equation*}

Our aim is to show that $B=A$, $($i.e. $Bf(z)=\frac{1}{2\pi i}\int_{\partial
\Delta }\frac{1}{w-z}H(z,w)f(w)dw).$ Fix $z\in \Delta $ and $f\in
L^{2}(\partial \Delta )$. Then%
\begin{equation*}
\left\Vert Bf-\sum_{k=0}^{m}\sum_{l=0}^{n}a_{kl}M_{k}C_{\Delta
}N_{l}f\right\Vert \rightarrow 0
\end{equation*}

i.e.%
\begin{equation*}
Bf=\lim_{m,n\rightarrow \infty
}\sum_{k=0}^{m}\sum_{l=0}^{n}a_{kl}M_{k}C_{\Delta }N_{l}f.
\end{equation*}%
Hence 
\begin{equation*}
\sum_{k=0}^{m}\sum_{l=0}^{n}a_{kl}M_{k}C_{\Delta }N_{l}f(z)\rightarrow Bf(z)
\end{equation*}%
Now%
\begin{eqnarray*}
Bf(z) &=&\lim_{m,n\rightarrow \infty
}\sum_{k=0}^{m}\sum_{l=0}^{n}a_{kl}M_{k}C_{\Delta }N_{l}f(z) \\
&=&\lim_{m,n\rightarrow \infty }\frac{1}{2\pi i}\sum_{k=0}^{m}\sum_{l=0}^{n}%
\int_{\partial \Delta }\frac{a_{kl}z^{k}w^{l}f(w)}{w-z}dw \\
&=&\lim_{m,n\rightarrow \infty }\frac{1}{2\pi i}\int_{\partial \Delta
}\sum_{k=0}^{m}\sum_{l=0}^{n}\frac{a_{kl}z^{k}w^{l}f(w)}{w-z}dw \\
&=&\frac{1}{2\pi i}\int_{\partial \Delta }\frac{f(w)}{w-z}\sum_{0}^{\infty
}\sum_{0}^{\infty }a_{mn}z^{m}w^{n}dw\qquad \text{(from equation \ref%
{equ:yyttr})} \\
&=&\frac{1}{2\pi i}\int_{\partial \Delta }\frac{f(w)}{w-z}H(z,w)dw \\
&=&Af(z).
\end{eqnarray*}%
So for $z\in \Delta $ and $f\in L^{2}(\partial \Delta )$, $Bf(z)=Af(z)$.
Hence $B=A$. Therefore, since $B$ is a continuous operator it follows that $%
A=U_{\psi }C_{D}\widetilde{U}_{\varphi }$ is a continuous operator.
Therefore $C_{D}$ is a continuous operator.


\begin{thebibliography}{99}
\bibitem{Bell} \label{bib:Bell}Bell, S.R., The Cauchy Transform, Potential
Theory and Conformal Mapping, Chapman and Hall/CRC, (2nd Edition, 2015)

\bibitem{Cima} \label{bib:Cima}Cima, J.A., Matheson, A.L., Ross W.T., The
Cauchy Transform, American Mathematical Society (2006).

\bibitem{ Coifman} \label{Coifman} Coifman, R.R., Jones, P.W., Semmes, S.,
Two elementary proofs of the L\symbol{94}2 boundedness of Cauchy integrals
on Lipschitz curves. J. Amer. Math.Soc. 2 (1989), no.3, 353-364.

\bibitem{dur} \label{bib:dur}Duren~P.L., Theory of $H^{p}$\ Spaces,
(Academic Press, 1970).

\bibitem{dur1} \label{bib:dur1}Duren~P.L., Smirnov Domains, Journal of
Mathematical Sciences, Volume 63, Number 2 / January, (1993), 167-170.

\bibitem{Dyn} \label{Dyn}Dyn'kin, E. M. Methods of the Theory of singular
integrals: Littlewood-Paley theory and its applications, Commutative
harmonic analysis IV, Volume 42 of the series Encyclopaedia of Mathematical
Sciences pp 97-194, Springer, Berlin, (1992).

\bibitem{Edmunds} \label{bib:Edmunds}Edmunds, D.E., Kokilashvili, V.,
Meskhi, A., Bounded and Compact Integral Operators (Mathematics and Its
Applications), Springer, (2002).

\bibitem{go} \label{bib:go} Goluzin, G.M., Functions of a Complex Variable%
\textit{,} Amer. Math. Soc., Providence, RI, (translated from Russian, 1969).

\bibitem{Halmos} \label{bib:Halmos}Halmos, P.R., Sunder, V.S., Bounded
Integral Operators on $L^{2}$\ Spaces, Springer, (1978).

\bibitem{kh1} \label{bib:kh1}{Khavinson}, D., Factorization theorems for
certain classes of analytic functions in multiply connected domains, \textit{%
Pacific. J. Math}. 108 (1983) 295-318.

\bibitem{li2} \label{bib:li2} {G. Little}, Equivalences of positive integral
operators with rational kernels, Proc. London. Math. Soc.\emph{\ }(3) 62
(1991) 403-426.

\bibitem{Murai} \label{bib:Murai}Murai, T., A Real Variable Method for the
Cauchy Transform, and Analytic Capacity, Springer, (2008).

\bibitem{Okikiolu} \label{bib:Okikiolu}Okikiolu, G.O., Aspects of Bounded
Integral Operators in Lp Spaces, Academic Press Inc (1971).

\bibitem{Tolsa1} \label{bib:Tolsa1}Tolsa, X., Analytic Capacity, the Cauchy
Transform, and Non-homogeneous Calder\'{o}n-Zygmund Theory, Birkh\"{a}user,
(2016).

\bibitem{Tolsa2} \label{bib:Tolsa2}Tolsa, X., Rectifiable Measures, Square
Functions Involving Densities, and the Cauchy Transform, American
Mathematical Society, (2017).
\end{thebibliography}
\end{document}